\documentclass[12pt, twoside, leqno]{article}
\usepackage{amsmath,amsthm}
\usepackage{amssymb,latexsym}
\usepackage{enumerate}

\newtheorem{thm}{Theorem}[section]

\newtheorem{cor}[thm]{Corollary}
\newtheorem{lem}[thm]{Lemma}

\theoremstyle{definition}

\newtheorem{rem}[thm]{Remark}

\numberwithin{equation}{section}

\usepackage{mathrsfs}
\usepackage{CJK}
\usepackage{amsmath}
\usepackage{fancyhdr}

\begin{document}

\baselineskip=17pt

\title{Polynomial maps with nilpotent Jacobians in dimension three I}
\author{ Dan Yan \footnote{ The author are supported by the Natural Science Foundation of Hunan Province (Grant No.2016JJ3085), the National Natural Science Foundation of China (Grant No.11601146) and the Construct Program of the Key Discipline in Hunan Province.}\\
Key Laboratory of HPCSIP,\\ College of Mathematics and Computer Science,\\
 Hunan Normal University, Changsha 410081, China \\
\emph{E-mail:} yan-dan-hi@163.com \\
}
\date{}

\maketitle

\renewcommand{\thefootnote}{}

\renewcommand{\thefootnote}{\arabic{footnote}}
\setcounter{footnote}{0}

\begin{abstract} In this paper, we first prove that $u,v,h$ are linearly
dependent over ${\bf K}$ if $JH$ is nilpotent and $H$ has the form:
$H=(u(x,y,z),v(u,h),h(x,y))$ with $H(0)=0$ or
$H=(u(x,y),v(u,h),h(x,y,z))$ with $H(0)=0$. Then we classify
polynomial maps of the form $H=(u(x,y),v(x,y,z), h(x,y))$ in the
case that $JH$ is nilpotent and $(\deg_yu,\deg_yh)\leq 2$.

\end{abstract}
{\bf Keywords.} Jacobian Conjecture, Nilpotent Jacobian matrix, Polynomial maps\\
{\bf MSC(2010).} Primary 14E05;  Secondary 14A05;14R15 \vskip 2.5mm

\section{Introduction}

Throughout this paper, we will write ${\bf K}$ for algebraically
closed field and ${\bf K}[X]={\bf K}[x_1,x_2,\ldots,x_n]$ (${\bf
K}[\bar{X}]={\bf K}[x,y,z]$)for the polynomial algebra over ${\bf
K}$ with $n$ ($3$) indeterminates. Let
$F=(F_1,F_2,\ldots,F_n):{\bf{K}}^n\rightarrow{\bf{K}}^n$ be a
polynomial map, that is, $F_i\in{\bf{K}}[X]$ for all $1\leq i\leq
n$. Let $JF=(\frac{\partial F_i}{\partial x_j})_{n\times n}$ be the
Jacobian matrix of $F$.

The Jacobian Conjecture (JC) raised by O.H. Keller in 1939 in
\cite{1} states that a polynomial map
$F:{\bf{K}}^n\rightarrow{\bf{K}}^n$ is invertible if the Jacobian
determinant $\det JF$ is a nonzero constant. This conjecture has
been attacked by many people from various research fields, but it is
still open, even for $n\geq 2$. Only the case $n=1$ is obvious. For
more information about the wonderful 70-year history, see \cite{2},
\cite{3}, and the references therein.

In 1980, S.S.S.Wang (\cite{4}) showed that the JC holds for all
polynomial maps of degree 2 in all dimensions (up to an affine
transformation). The most powerful result is the reduction to degree
3, due to H.Bass, E.Connell and D.Wright (\cite{2}) in 1982 and
A.Yagzhev (\cite{5}) in 1980, which asserts that the JC is true if
the JC holds for all polynomial maps $X+H$, where $H$ is homogeneous
of degree 3. Thus, many authors study these maps and led to pose the
following problem.

 {\em (Homogeneous) dependence problem.} Let $H=(H_1,\ldots,H_n)\in
{\bf K}[X]$ be a (homogeneous) polynomial map of degree $d$ such
that $JH$ is nilpotent and $H(0)=0$. Whether $H_1,\ldots,H_n$ are
linearly dependent over ${\bf K}$?

The answer to the above problem is affirmative if rank$JH\leq 1$
(\cite{2}). In particular, this implies that the dependence problem
has an affirmative answer in the case $n=2$. D. Wright give an
affirmative answer when $H$ is homogeneous of degree 3 in the case
$n=3$ (\cite{6}) and the case $n=4$ is solved by Hubbers in
\cite{7}. M. de Bondt and A. van den Essen give an affirmative
answer to the above problem in the case $H$ is homogeneous and $n=3$
(\cite{8}). A. van den Essen finds the first counterexample in
dimension three for the dependence problem (\cite{9}). M. de Bondt
give a negative answer to the homogeneous dependence problem for
$d\geq 3$. In particular, he constructed counterexamples to the
problem for all dimensions $n\geq 5$ (\cite{10}). In \cite{18}, M.
Chamberland and A. van den Essen classify all polynomial maps of the
form $H=(u(x,y),v(x,y,z),h(u(x,y),v(x,y,z)))$ with $JH$ nilpotent.
In particular, they show that all maps of this form with $H(0)=0$,
$JH$ nilpotent and $H_1, H_2, H_3$ are linearly independent has the
same form as the counterexample that gave by A. van den Essen in
\cite{9} (up to a linear coordinate change). We classify all
polynomial maps of the form $H=(u(x,y),v(x,y,z),h(x,y,z))$ in the
case that $JH$ is nilpotent and $\deg_zv\leq 3$,
$(\deg_yu(x,y),\deg_yh(x,y,z))=1$ (\cite{13}).

In section 2, we prove that $u,v,h$ are linearly dependent over
${\bf K}$ if $JH$ is nilpotent and $H$ has the form:
$H=(u(x,y,z),v(u,h),h(x,y))$ with $H(0)=0$ or
$H=(u(x,y),v(u,h),h(x,y,z))$ with $H(0)=0$. Then, in section 3, we
classify all polynomial maps of the form
$H=(u(x,y),v(x,y,z),h(x,y))$ in the case that $JH$ is nilpotent and
$(\deg_yu(x,y),\deg_yh(x,y))\leq 2$ or $\deg_yu(x,y)$ or
$\deg_yh(x,y)$ is a prime number. In Theorem 3.11, we prove that
$u,v,h$ are linearly dependent over ${\bf K}$ if $JH$ is nilpotent
and $H$ has the form: $H=(u(x,y),v(x,y,z),h(x,y))$ with $H(0)=0$ and
$u$ is homogeneous. The main results in the paper are Theorem 2.4,
Theorem 2.6, Theorem 3.3, Theorem 3.5 and Theorem 3.11. We define
$Q_{x_i}:=\frac{\partial Q}{\partial x_i}$ and that $\deg_y f$ is
the highest degree of $y$ in $f$.

\section{Polynomial maps of the form $H=(u(x,y,z),\allowbreak v(u,h), h(x,y))$}

In this section, we prove that $H$ are linearly dependent over ${\bf
K}$ if $JH$ is nilpotent and $H$ has the form:
$H=(u(x,y,z),v(u,h),h(x,y))$ with $H(0)=0$ or
$H=(u(x,y),v(u,h),h(x,y,z))$ with $H(0)=0$.

\begin{thm}
Let $H=(u(x,y,z), v(h(x,y)), h(x,y))$ be a polynomial map with
$H(0)=0$. If $JH$ is nilpotent, then $u,v,h$ are linearly dependent.
\end{thm}
\begin{proof}
If $\deg_zu=0$, then the conclusion follows from Proposition 2.1 in
\cite{13}. Suppose that $\deg_zu\geq 1$. Since $JH$ is nilpotent, we
have the following equations:
\begin{equation}
\nonumber
  \left\{ \begin{aligned}
  u_x+v'(h)h_y = 0~~~~~~~~~~~~~~~~~~~~~~~~~~~~~(2.1) \\
  u_xv'(h)h_y-u_yv'(h)h_x-u_zh_x = 0~~~~~~~(2.2) \\
                          \end{aligned} \right.
  \end{equation}
Let $u=u_mz^m+u_{m-1}z^{m-1}+\cdots+u_1z+u_0$ with $u_m\neq 0$. It
follows from equation $(2.1)$ that
$$u_{mx}z^m+u_{(m-1)x}z^{m-1}+\cdots+u_{1x}z+u_{0x}+v'(h)h_y=0$$
We always view that the polynomials are in ${\bf K}[x,y][z]$ with
coefficients in ${\bf K}[x,y]$. Comparing the coefficients of the
degree of $z$ of the above equation, we have
$u_{mx}=u_{(m-1)x}=\cdots=u_{1x}=0$ and
$u_{0x}=-v'(h)h_y=-\frac{\partial(v(h))}{\partial y}$. It follows
from equation $(2.2)$ that $v'(h)(u_xh_y-u_yh_x)=u_zh_x$. That is,
$$v'(h)[u_{0x}h_y-(u_{my}z^m+u_{(m-1)y}z^{m-1}+\cdots+u_{1y}z+u_{0y})h_x]$$
$$=h_x(mu_mz^{m-1}+(m-1)u_{m-1}z^{m-2}+\cdots+u_1)~~~~~~~~~(2.3)$$
Comparing the coefficients of $z^m$ of equation $(2.3)$, we have
$v'(h)u_{my}h_x=0$. Thus, we have $v'(h)=0$ or $h_x=0$ or
$u_{my}=0$.

Case (i) If $v'(h)=0$, then $v(h)=c$. Since $H(0)=0$, we have
$c=v(h(0,0))=v(0)=0$. Thus, $u,v,h$ are linearly dependent.

Case (ii) If $h_x=0$, then it follows from equation $(2.2)$ that
$u_xv'(h)h_y=0$. Thus, we have $u_x=0$ or $h_y=0$ or $v'(h)=0$.

$(1)$ If $v'(h)=0$, then it reduces to Case (i).

$(2)$ If $h_y=0$, then $h=0$ because $H(0)=0$. Thus, $u,v,h$ are
linearly dependent.

$(3)$ If $u_x=0$, then it follows from equation $(2.1)$
$v'(h)h_y=0$. That is, $v'(h)=0$ or $h_y=0$. Then it reduces to Case
(i) and Case (ii) $(2)$ respectively.

Case (iii) If $u_{my}=0$, then $u_m\in{\bf K}^*$.

Suppose $m\geq 2$. Then comparing the coefficients of $z^{m-1}$ of
equation $(2.3)$, we have $$-v'(h)u_{(m-1)y}h_x=mu_mh_x.$$ We can
assume that $h_x\neq0$, otherwise, it follows from Case (ii) that
$u,v,h$ are linearly dependent. Thus, we have
$-v'(h)u_{(m-1)y}=mu_m$. Therefore, we have $v'(h)\in {\bf K}^*$ and
$u_{(m-1)y}\in {\bf K}^*$. That is, $v(h)=c_1h+c_0$ for $c_1\in {\bf
K}^*$, $c_0\in {\bf K}$. Since $H(0)=0$, we have that $c_0=0$. Thus,
$u,v,h$ are linearly dependent.

Suppose $m=1$. Then equation $(2.3)$ has the following form:
$$v'(h)(u_{0x}h_y-u_{0y}h_x)=u_1h_x~~~~~~~~~~~~(2.4)$$
where $u_1\in {\bf K}^*$. Comparing the degree of $x$ of the above
equation, we have $\deg_hv(h)\leq 1$ in the case
$u_{0x}h_y-u_{0y}h_x\neq 0$. Then we have
$v(h)=\bar{c}h+\bar{\bar{c}}$ for $\bar{c},\bar{\bar{c}}\in{\bf K}$.
Since $H(0)=0$, we have $\bar{\bar{c}}=v(h(0,0))=v(0)=0$. That is,
$v(h)=\bar{c}h$. Thus, $u,v,h$ are linearly dependent. If
$u_{0x}h_y-u_{0y}h_x = 0$, then it follows from equation $(2.4)$
that $h_x=0$. Then it reduces to Case (ii).
\end{proof}

\begin{rem}
If $H(0)\neq 0$, then it easy to compute that there exists
$H=(u,v,h)$ such that $u, v, h$ are linearly independent and $JH$ is
nilpotent: $H=(u(y), c, h(x,y))$ with $c\neq 0$, $u(y)\in {\bf
K}[y]/{\bf K}$ and $h(x,y)\in {\bf K}[x,y]/{\bf K}[y]$.
\end{rem}

\begin{cor}
Let $H=(u(h), v(x,y,z), h(x,y))$ be a polynomial map with $H(0)=0$.
If $JH$ is nilpotent, then $u,v,h$ are linearly dependent.
\end{cor}
\begin{proof}
Let \[T=\left(
  \begin{array}{ccc}
    0 & 1 & 0 \\
    1 & 0 & 0 \\
    0 & 0 & 1 \\
  \end{array}
\right).\] Then $T^{-1}HT=(v(y,x,z),u(h(y,x)),h(y,x))$. Since $JH$
is nilpotent, we have that $J(T^{-1}HT)=T^{-1}JHT$ is nilpotent. It
follows from Theorem 2.1 that $u,v,h$ are linearly dependent.
\end{proof}

In the proof of the following theorem, our main goal is to reduce it
to Theorem 2.1. We divide the proof into two parts according to the
degree of $u$ in $v$.

\begin{thm}
Let $H=(u(x,y,z),v(u,h),h(x,y))$ be a polynomial map over ${\bf
K}[x,y,z]$ with $H(0)=0$. If $JH$ is nilpotent, then $u,v,h$ are
linearly dependent over ${\bf K}$.
\end{thm}
\begin{proof}
Since $JH$ is nilpotent, we have the following equations:
\begin{equation}
\nonumber
  \left\{ \begin{aligned}
  u_x+v_hh_y+v_uu_y = 0~~~~~~~~~~~~~~~~~~~~~~~~~~~~~(2.5) \\
  u_xh_yv_h-u_yh_xv_h-u_zh_x-h_yv_uu_z = 0~~~~~~~(2.6) \\
                          \end{aligned} \right.
  \end{equation}
  Let $u=u_m(x,y)z^m+u_{m-1}(x,y)z^{m-1}+\cdots+u_1(x,y)z+u_0(x,y)$,
  $v(u,h)=v_n(h)u^n+v_{n-1}(h)u^{n-1}+\cdots+v_1(h)u+v_0(h)$ with $u_mv_n\neq
  0$. If $m=0$, then it follows from Proposition 2.1 in \cite{13}
  that $u,v,h$ are linearly dependent. If $n=0$, then the conclusion
  follows from Theorem 2.1. Thus, we can assume that $mn\geq 1$. It
  follows from equation $(2.5)$ that
\\$u_{mx}z^m+u_{(m-1)x}z^{m-1}+\cdots+u_{1x}z+u_{0x}+[v_n'(h)u^n+v_{n-1}'(h)u^{n-1}+\cdots+v_1'(h)u+v_0'(h)]h_y
  +(u_{my}z^m+u_{(m-1)y}z^{m-1}+\cdots+u_{1y}z+u_{0y})(nv_n(h)u^{n-1}+(n-1)v_{n-1}(h)u^{n-2}+\cdots+v_1(h))=0~~~~~~~~~~~~~~~~~~~~~~~~~~~~~~~~~~~~~~~~~~~~~~~~~(2.7)$
\\We always view that the polynomials are in ${\bf K}[x,y][z]$ with
coefficients in ${\bf K}[x,y]$ in the following arguments.\\

If $n\geq 2$, then comparing the coefficients of $z^{mn}$ of
equation $(2.7)$, we have the following equation:
$$v_n'(h)h_yu_m+nv_n(h)u_{my}=0$$
That is,
$$\frac{v_n'(h)h_y}{v_n(h)}=-n\frac{u_{my}}{u_m}~~~~~~~~~~~~~~~~~~~(2.8)$$
Suppose $u_{my}\neq 0$. Then $v_n'(h)h_y\neq 0$. Thus, we have
$v_n(h)u_m^n=e^{c(x)}$ by integrating the two sides of equation
$(2.8)$ with respect to $y$, where $c(x)$ is a function of $x$.
Since $v_n(h), u_m\in {\bf K}[x,y]$ and $e^{c(x)}$ is a function of
$x$, so we have $u_m, v_n(h)\in {\bf K}[x]$. This is a
contradiction! Therefore, we have $u_{my}=0$ and $v_n'(h)h_y=0$.
That is, $u_{my}=0$ and $v_n'(h)=0$
or $u_{my}=0$ and $h_y=0$.\\

Case I If $u_{my}=0$ and $h_y=0$, then we have
$u_{(m-1)y}=\cdots=u_{1y}=0$ by comparing the coefficients of
$z^{m(n-1)+i}$ for $i=m-1, \ldots, 2, 1$ of equation $(2.7)$. It
follows from equation $(2.6)$ that
$[u_{0y}(v_n'(h)u^n+v_{n-1}'(h)u^{n-1}+\cdots+v_1'(h)u+v_0'(h))+mu_mz^{m-1}+(m-1)u_{m-1}z^{m-2}+\cdots+u_1]h_x=0$
Thus, we have $h_x=0$ or
$u_{0y}(v_n'(h)u^n+v_{n-1}'(h)u^{n-1}+\cdots+v_1'(h)u+v_0'(h))+mu_mz^{m-1}+(m-1)u_{m-1}z^{m-2}+\cdots+u_1=0~~~~~~~~~~~~~~~~~~~~~~~~~~~~~~~~~~~~~~~~~~~~~~~~~~~~~~~~~~~~(2.9)$

(i) If $h_x=0$, then $h=0$ because $H(0)=0$. Thus, $u,v,h$ are
linearly dependent.

(ii) If equation $(2.9)$ is true, then $u_{0y}v_n'(h)=0$ by
comparing the coefficients of $z^{nm}$ of equation $(2.9)$. Thus, we
have $u_{0y}=0$ or $v_n'(h)=0$. If $u_{0y}=0$, then it follows from
equation $(2.9)$ that $u_m=0$. This is a contradiction! If
$v_n'(h)=0$, then it reduces to the following case.\\

Case II If $u_{my}=0$ and $v_n'(h)=0$, then we have
$u_{(m-1)y}=\cdots=u_{1y}=0$ by comparing the coefficients of
$z^{m(n-1)+i}$ for $i= m-1, \ldots, 2,1$ of equation $(2.7)$.
Comparing the coefficients of $z^{im}$ for $i= n-1, n-2,\ldots,2$ of
equation $(2.7)$, we have
$$v_i'(h)h_y+(i+1)v_{i+1}(h)u_{0y}=0~~~~~~~~~~~~~~~(2.10)$$
and we have
$$u_{jx}+u_j(v_1'(h)h_y+2v_2(h)u_{0y})=0~~~~~~~~~~~~~~(2.11)$$
by comparing the coefficients of $z^j$ for $j=m,m-1,\ldots,1$. Then
equation $(2.7)$ has the following form:
$$u_{0x}+u_0(v_1'(h)h_y+2v_2(h)u_{0y})+v_0'(h)h_y+u_{0y}v_1(h)=0~~~~~~~~~~(2.12)$$
It follows from equation $(2.6)$ that
$(v_{n-1}'(h)u^{n-1}+v_{n-2}'(h)u^{n-2}+\cdots+v_1'(h)u+v_0'(h))[(u_{mx}z^m+u_{(m-1)x}z^{m-1}+\cdots+u_{1x}z+u_{0x})h_y-u_{0y}h_x]-h_x(mu_mz^{m-1}+
(m-1)u_{m-1}z^{m-2}+\cdots+u_1)-h_y(nv_n(h)u^{n-1}+(n-1)v_{n-1}(h)u^{n-2}+\cdots+v_1(h))(mu_mz^{m-1}+(m-1)u_{m-1}z^{m-2}+\cdots+u_1)=0
~~~~~~~~~~~~~~~~~~~~~~(2.13)$ \\
Comparing the coefficients of $z^{nm}$ of equation $(2.13)$, we have
$v_{n-1}'(h)u_{mx}=0$. Thus, we have $v_{n-1}'(h)=0$ or
$u_{mx}=0$.\\

If $v_{n-1}'(h)=0$ and $m\geq 2$, then we have $h_ynv_n(h)mu_m^n=0$
by comparing the coefficients of $z^{(n-1)m+m-1}$ of equation
$(2.13)$. Thus, we have $h_y=0$. It follows from equation $(2.10)$
($i=m-1$) that $u_{0y}=0$. Comparing the coefficients of $z^{m-1}$
of equation $(2.13)$, we have $h_x=0$. Since $H(0)=0$, so we have
$h=0$. Thus, $u,v,h$ are linearly dependent.\\

If $v_{n-1}'(h)=0$ and $m = 1$, then we have
$$h_y(v_{n-2}'(h)u_{1x}-nv_n(h)u_1^2)=0$$
by comparing the coefficients of $z^{n-1}$ of equation $(2.13)$.
Thus, we have $h_y=0$ or $v_{n-2}'(h)u_{1x}-nv_n(h)u_1^2=0$.

(a) If $h_y=0$, then we have $u_{0y}=0$ by equation $(2.10)$
($i=m-1$). It follows from equation $(2.13)$ that $h_x=0$. Since
$H(0)=0$, we have $h=0$. Thus, $u,v,h$ are linearly dependent.

(b) If $v_{n-2}'(h)u_{1x}-nv_n(h)u_1^2=0$, then we have
$$v_{n-2}'(h)u_{1x}=nv_nu_1^2~~~~~~~~~~~~~~~~~~(2.14)$$
Since $u_{1y}=0$, we have $u_1\in{\bf K}[x]$. Thus, $nv_nu_1^2
\in{\bf K}[x]$. It follows from equation $(2.14)$ that
$v_{n-2}'(h)\in {\bf K}$ or $h_y=0$.\\
If $h_y=0$, then it reduces to (a). If $v_{n-2}'(h)\in {\bf K}$,
then we have $v_nu_1=0$ by comparing the degree of $x$ of two sides
of equation $(2.14)$. This is a contradiction!\\

If $u_{mx}=0$, then we have $v_n(h), u_m\in {\bf K}^*$. It follows
from equation $(2.11)$ ($j=m$) that
$$v_1'(h)h_y+2v_2(h)u_{0y}=0~~~~~~~~~~~~~(2.15)$$
Thus, we have $u_{jx}=0$ for $1\leq j\leq m-1$ by substituting
equation $(2.15)$ to equation $(2.11)$ for $1\leq j\leq m-1$.
Therefore, equation $(2.12)$ has the following form:
$$u_{0x}+v_0'(h)h_y+u_{0y}v_1(h)=0~~~~~~~~~~~(2.16)$$
and equation $(2.13)$ has the following form:
$(v_{n-1}'(h)u^{n-1}+v_{n-2}'(h)u^{n-2}+\cdots+v_1'(h)u+v_0'(h))(u_{0x}h_y-u_{0y}h_x)-h_x(mu_mz^{m-1}+
(m-1)u_{m-1}z^{m-2}+\cdots+u_1)-h_y(nv_n(h)u^{n-1}+(n-1)v_{n-1}(h)u^{n-2}+\cdots+v_1(h))(mu_mz^{m-1}+(m-1)u_{m-1}z^{m-2}+\cdots+u_1)=0
~~~~~~~~~~~~~~~~~~~~~~~~~~~~~~~~~~~~~~~~~~~~~~~~~~~~~~~~~~(2.17)$

If $m\geq 2$, then we have $h_ynv_n(h)mu_m^n=0$ by comparing the
coefficients of $z^{m(n-1)+m-1}$ of equation $(2.17)$. Thus, we have
$h_y=0$. Then, with the same arguments as in (a), we have that
$u,v,h$ are linearly dependent.

If $m=1$, then $u=u_1z+u_0(x,y)$,
$v(u,h)=v_n(h)u^n+\cdots+v_1(h)u+v_0(h)$ with $u_1, v_n(h)\in {\bf
K}^*$. Thus, equation $(2.17)$ has the following form:
$(v_{n-1}'(h)u^{n-1}+v_{n-2}'(h)u^{n-2}+\cdots+v_1'(h)u+v_0'(h))(u_{0x}h_y-u_{0y}h_x)-h_xu_1-h_yu_1(nv_n(h)u^{n-1}+(n-1)v_{n-1}(h)u^{n-2}+\cdots+v_1(h))=0
~~~~~~~~~~~~~~~~~~~~~~~~~~(2.18)$\\
Comparing the coefficients of $z^{im}$ of equation $(2.18)$ for
$i=n-1,\ldots,1$, we have the following equations:
$$v_i'(h)(u_{0x}h_y-u_{0y}h_x)-(i+1)v_{i+1}(h)u_1h_y=0~~~~~~(2.19)$$
Then equation $(2.18)$ has the following form:
$$v_0'(h)(u_{0x}h_y-u_{0y}h_x)-h_xu_1-v_1(h)u_1h_y=0~~~~~~~~(2.20)$$
If $u_{0x}h_y-u_{0y}h_x=0$, then it follows from equation $(2.19)$
for $i=n-1$ that $h_y=0$. It follows from equation $(2.18)$ that
$h_x=0$. Since $H(0)=0$, we have $h=0$. Thus, $u,v,h$ are linearly
dependent.\\
If $u_{0x}h_y-u_{0y}h_x\neq 0$, then we have $v_{n-1}'(h)\in {\bf
K}$ by comparing the degree of $y$ of equation $(2.19)$ for $i=n-1$.
It follows from equation $(2.10)$ ($i=n-1$) that
$$u_0=-\frac{v_{n-1}(h)}{nv_n}+c_1(x)~~~~~~~~~~~~(2.21)$$
for some $c_1(x)\in {\bf K}[x]$. Substituting equation $(2.21)$ to
equation $(2.19)$ for $i=n-1$, we have
$v_{n-1}'(h)c_1'(x)h_y=nv_nu_1h_y$. If $h_y=0$, then it follows from
the arguments of (a) that $h=0$. Thus, $u,v,h$ are linearly
dependent. If $h_y\neq 0$, then $v_{n-1}'(h)c_1'(x)=nv_nu_1$. Since
$v_n,u_1\in{\bf K}^*$, so we have $v_{n-1}'(h), c_1'(x)\in {\bf
K}^*$. Thus, $c_1(x)=ax+b$ with $a\in{\bf K}^*$. Then we have
$$a=\frac{nv_nu_1}{v_{n-1}'(h)}~~~~~~~~~~~~~~~~~~~~~~~(2.22)$$
Thus, equation $(2.19)$ $(i=n-1)$ has the following form:
$$u_{0x}h_y-u_{0y}h_x=ah_y.$$
Substituting the above equation to equation $(2.20)$, we have the
following equation:
$$v_0'(h)ah_y=u_1h_x+v_1(h)u_1h_y$$
Substituting equation $(2.22)$ to the above equation, we have
$$v_0'(h)h_y=\frac{v_{n-1}'(h)}{nv_n}h_x+\frac{v_1(h)v_{n-1}'(h)}{nv_n}h_y~~~~~~~~(2.23)$$
Then substituting equation $(2.21)$ to equation $(2.16)$, we have
$$v_0'(h)h_y=\frac{v_{n-1}'(h)}{nv_n}h_x-a+\frac{v_1(h)v_{n-1}'(h)}{nv_n}h_y$$
Substituting equation $(2.23)$ to the above equation, we have $a=0$.
This is a contradiction!\\

If $n=1$, then
$u=u_m(x,y)z^m+u_{m-1}(x,y)z^{m-1}+\cdots+u_1(x,y)z+u_0(x,y)$,
$v(u,h)=v_1(h)u+v_0(h)$ with $u_mv_1\neq 0$. It follows from
equation $(2.6)$ that
$(v_1'(h)u+v_0'(h))[(u_{mx}z^m+u_{(m-1)x}z^{m-1}+\cdots+u_{1x}z+u_{0x})h_y-(u_{my}z^m+u_{(m-1)y}z^{m-1}+\cdots+u_{1y}z+u_{0y})h_x]-h_x(mu_mz^{m-1}+
(m-1)u_{m-1}z^{m-2}+\cdots+u_1)-h_yv_1(h)(mu_mz^{m-1}+(m-1)u_{m-1}z^{m-2}+\cdots+u_1)=0
~~~~~~~~~~~~~~~~~~~~~~~~~(2.24)$\\
Comparing the coefficients of $z^{2m}$ of equation $(2.24)$, we have
$v_1'(h)u_m(u_{mx}h_y-u_{my}h_x)=0$. Thus, we have $v_1'(h)=0$ or
$u_{mx}h_y-u_{my}h_x=0$.

If $v_1'(h)=0$, then $v_1(h)\in {\bf K}^*$. Let
 \[T_1=\left(
  \begin{array}{ccc}
    1 & 0 & 0 \\
    v_1 & 1 & 0 \\
    0 & 0 & 1 \\
  \end{array}
\right).\] Then
$T_1^{-1}HT_1=(u(x,y+v_1x),v_0(h(x,y+v_1x)),h(x,y+v_1x)):=(\bar{u},\bar{v},\bar{h})$.
It follows from Theorem 2.1 that $\bar{u},\bar{v},\bar{h}$ are
linearly dependent. Thus, $u,v,h$ are linearly dependent.

If $v_1'(h)\neq 0$, then $u_{mx}h_y-u_{my}h_x=0$ and
$u_{ix}h_y-u_{iy}h_x=0$ by comparing the coefficients of $z^{m+i}$
of equation $(2.24)$ for $i=m-1,m-2,\ldots,2,1,0$. That is,
$$u_xh_y-u_yh_x=0~~~~~~~~~~~~~~~~~~~~~~~~~(2.25)$$
Then equation $(2.24)$ has the following form: $(h_x+v_1h_y)u_z=0$.
Since $u_z\neq 0$, so we have $h_x+v_1h_y=0$. That is,
$$h_x=-v_1(h)h_y~~~~~~~~~~~~~~~~~~~~~~~~~~~~(2.26)$$
Substituting equation $(2.26)$ to equation $(2.25)$, we have
$(u_x+v_1(h)u_y)h_y=0$. Thus, we have $h_y=0$ or
$u_x+v_1(h)u_y=0$.\\
If $h_y=0$, then it follows from equation $(2.26)$ that $h_x=0$.
Since $H(0)=0$, we have $h=0$. Thus, $u,v,h$ are linearly
dependent.\\
If $u_x+v_1(h)u_y=0$, then it follows from equation $(2.5)$ that
$$u_x+v_hh_y+v_1(h)u_y=0$$
Thus, we have $v_hh_y=0$. Since $v_h\neq 0$, so we have $h_y=0$. It
follows from the above arguments that $u,v,h$ are linearly
dependent.
\end{proof}

\begin{cor}
Let $H=(u(v,h),v(x,y,z),h(x,y))$ be a polynomial map with $H(0)=0$.
If $JH$ is nilpotent, then $u,v,h$ are linearly dependent.
\end{cor}
\begin{proof}
Let
\[T=\left(
  \begin{array}{ccc}
    0 & 1 & 0 \\
    1 & 0 & 0 \\
    0 & 0 & 1 \\
  \end{array}
\right).\] Then $T^{-1}HT=(v(y,x,z),u(v,h),h(y,x))$. Since $JH$ is
nilpotent, we have that $J(T^{-1}HT)$ is nilpotent. It follows from
Theorem 2.2 that $u,v,h$ are linearly dependent.
\end{proof}

\begin{thm}
Let $H=(u(x,y),v(u,h),h(x,y,z))$ be a polynomial map with $H(0)=0$.
If $JH$ is nilpotent, then $u,v,h$ are linearly dependent over ${\bf
K}$.
\end{thm}
\begin{proof}
Since $JH$ is nilpotent, we have the following equations:
\begin{equation}
\nonumber
  \left\{ \begin{aligned}
  u_x+v_uu_y+v_hh_y+h_z = 0~~~~~~~~~~~~~~~~~~~~~~(2.27) \\
  u_xv_hh_y-u_yv_hh_x+h_zv_uu_y+u_xh_z = 0~~~~~~~(2.28) \\
                          \end{aligned} \right.
  \end{equation}
Let $h(x,y,z)=h_dz^d+h_{d-1}z^{d-1}+\cdots+h_1z+h_0$,
$v(u,h)=v_n(u)h^n+v_{n-1}(u)h^{n-1}+\cdots+v_1(u)h+v_0(u)$ with
$h_dv_n(u)\neq 0$.

If $nd=0$, then it follows from Proposition 2.1 that $u,v,h$ are
linearly dependent. Thus, we can assume that $n\geq 1$ and $d\geq
1$. It follows from equation $(2.27)$ that
$u_x+(v_n'(u)h^n+v_{n-1}'(u)h^{n-1}+\cdots+v_1'(u)h+v_0'(u))u_y+
(nv_n(u)h^{n-1}+(n-1)v_{n-1}(u)h^{n-2}+\cdots+v_1(u))(h_{dy}z^d+h_{(d-1)y}z^{d-1}+\cdots+h_{1y}z+h_{0y})
+dh_dz^{d-1}+(d-1)h_{d-1}z^{d-2}+\cdots+h_1=0~~~~~~~~~~~~~~~~~~~~~~~~~~~~~~~~~~~~~~~~(2.29)$\\
We always view that the polynomials are in ${\bf K}[x,y][z]$ with
coefficients in ${\bf K}[x,y]$ in the following arguments. Comparing
the coefficients of $z^{dn}$ of equation $(2.29)$, we have
$$v_n'(u)h_d^nu_y+nv_n(u)h_d^{n-1}h_{dy}=0.$$
That is,
$$\frac{v_n'(u)u_y}{v_n(u)}=-n\frac{h_{dy}}{h_d}~~~~~~~~~~~~~~~~~~~~~(2.30)$$
Suppose $h_{dy}\neq 0$. Then $v_n'(u)u_y\neq 0$. Thus, we have
$v_n(u)h_d^n=e^{c(x)}$ by integrating the two sides of equation
$(2.30)$ with respect to $y$, where $c(x)$ is a function of $x$.
Since $v_n(u),h_d\in {\bf K}[x,y]$ and $e^{c(x)}$ is a function of
$x$, we have $v_n(u),h_d\in{\bf K}[x]$. This is a contradiction!
Thus, we have that $h_{dy}=0$ and $v_n(u)\in {\bf K}^*$ or
$h_{dy}=0$ and
$u_y=0$.\\

If $h_{dy}=0$ and $u_y=0$, then equation $(2.27)$ has the following
form: $$u_x+v_hh_y+h_z=0~~~~~~~~~~~~~~~~~~(2.31)$$ It follows from
equation $(2.28)$ that $u_x(v_hh_y+h_z)=0$. Thus, we have $u_x=0$ or
$v_hh_y+h_z=0$.

(A) If $u_x=0$, then $u=0$ because $H(0)=0$. Thus, $u,v,h$ are
linearly dependent.

(B) If $v_hh_y+h_z=0$, then it follows from equation $(2.31)$ that
$u_x=0$. It reduces to (A).\\

If $h_{dy}=0$ and  $v_n(u)\in {\bf K}^*$, then suppose that $n\geq
2$, we have $h_{iy}=0$ for $i=d-1,\ldots,1$ by comparing the
coefficients of $z^{d(n-1)+i}$ of equation $(2.29)$. Thus, equation
$(2.29)$ has the following
form:\\
 $u_x+(v_{n-1}'(u)h^{n-1}+\cdots+v_1'(u)h+v_0'(u))u_y+
(nv_n(u)h^{n-1}+(n-1)v_{n-1}(u)h^{n-2}+\cdots+v_1(u))h_{0y}
+dh_dz^{d-1}+(d-1)h_{d-1}z^{d-2}+\cdots+h_1=0~~~~~~~~~(2.32)$\\
Comparing the coefficients of $z^{dj}$ for $j=n-1, n-2,\ldots,1$, we
have the following equations:
$$v_j'(u)u_y+(j+1)v_{j+1}(u)h_{0y}=0~~~~~~~~~~~~~~~~(2.33)$$
Thus, we have $d=1$ by comparing the coefficients of $z^{d-1}$. Then
equation $(2.32)$ has the following form:
$$v_0'(u)u_y=-u_x-v_1(u)h_{0y}-h_1~~~~~~~~~~~~~~~~~~(2.34)$$
Consequently, we have $h(x,y,z)=h_1z+h_0$,
$v(u,h)=v_n(u)h^n+v_{n-1}(u)h^{n-1}+\cdots+v_1(u)h+v_0(u)$ with
$v_nh_1\neq 0$ and $h_{1y}=0$, $v_n(u)\in {\bf K}^*$. It follows
from equation $(2.28)$ that
$u_x(nv_nh^{n-1}+(n-1)v_{n-1}(u)h^{n-2}+\cdots+v_1(u))h_{0y}-u_y(nv_nh^{n-1}+(n-1)v_{n-1}(u)h^{n-2}+\cdots+v_1(u))
(h_{1x}z+h_{0x})+h_1(v_{n-1}'(u)h^{n-1}+\cdots+v_1'(u)h+v_0'(u))u_y+u_xh_1=0
~~~~~~~~~~~~~~~~~~~~~~~~~~~~~~~~~~~~~~~~~~~~~~~~~~~~~(2.35)$\\
Then we have $u_ynv_nh_1^{n-1}h_{1x}=0$ by comparing the
coefficients of $z^n$ of equation $(2.35)$. Thus, we have $u_y=0$ or
$h_{1x}=0$.

If $u_y=0$, then it reduces to the former case.

If $h_{1x}=0$, then $h_1\in{\bf K}^*$. Comparing the coefficients of
$z^j$ for $j=n-1,\ldots,1$ of equation $(2.35)$, we have
$$(j+1)v_{j+1}(u)(u_xh_{0y}-u_yh_{0x})+h_1v_j'(u)u_y=0~~~~~~~~~~(2.36)$$
Then equation $(2.35)$ has the following form:
$$v_1(u)(u_xh_{0y}-u_yh_{0x})+h_1v_0'(u)u_y+u_xh_1=0~~~~~~~~~~~~(2.37)$$
Substituting equations $(2.33)$ to equation $(2.36)$ for $j=n-1$, we
have the follows equation:
$$u_xh_{0y}-u_yh_{0x}=h_1h_{0y}~~~~~~~~~~~~~~~~~~~~~~~~~~~~~~~~~~(2.38)$$
Substituting equations $(2.34)$, $(2.38)$ to equation $(2.37)$, we
have
$$v_1(u)h_1h_{0y}+h_1(-u_x-v_1(u)h_{0y}-h_1)+u_xh_1=0$$ That is,
$h_1=0$. This is a contradiction! Therefore, $u,v,h$ are
linearly dependent.\\

If $h_{dy}=0$ and  $v_n(u)\in {\bf K}^*$ and $n=1$, then
$v(u,h)=v_1(u)h+v_0(h)$ and $v_1:=v_1(u)\in {\bf K}^*$. Then
equation $(2.29)$ has the following form:
$u_x+v_0'(u)u_y+v_1(h_{(d-1)y}z^{d-1}+\cdots+h_{1y}z+h_{0y})+dh_dz^{d-1}+(d-1)h_{d-1}z^{d-2}+\cdots+h_1=0$~~~~~~~~~~(2.39)
Then we have
$$v_1h_{(d-1)y}+dh_d=0~~~~~~~~~~~~~~~~~~~~~~~~~~~~~~(2.40)$$
by comparing the coefficients of $z^{d-1}$ of equation $(2.39)$. It
follows from equation $(2.28)$ that\\
$u_xv_1(h_{(d-1)y}z^{d-1}+\cdots+h_{1y}z+h_{0y})-u_yv_1(h_{dx}z^d+h_{(d-1)x}z^{d-1}+\cdots+h_{1x}z+h_{0x})+
(dh_dz^{d-1}+(d-1)h_{d-1}z^{d-2}+\cdots+h_1)(v_0'(u)u_y+u_x)=0$~~~~~~~~~~~~~~~~~~~~~(2.41)
Comparing the coefficient of $z^d$ of equation $(2.41)$, we have
$u_yh_{dx}=0$. That is, $u_y=0$ or $h_{dx}=0$.

(I) If $u_y=0$, then $u_x=0$ because $JH$ is nilpotent. Since
$u(0,0)=0$, we have $u=0$. Thus, $u,v,h$ are linearly dependent.

(II) If $h_{dx}=0$, then we have $h_d\in{\bf K}^*$ and
$$v_1(u_xh_{(d-1)y}-u_yh_{(d-1)x})+dh_d(u_x+v_0'(u)u_y)=0~~~~~~~~~(2.42)$$
by comparing the coefficients of $z^{d-1}$ of equation $(2.41)$.
Substituting equation $(2.40)$ to equation $(2.42)$, we have
$$(dh_dv_0'(u)-v_1h_{(d-1)x})u_y=0~~~~~~~~~~~~~~~~~(2.43)$$
It follows from equation $(2.40)$ that $h_{d-1}=-dv_1^{-1}h_d\cdot
y+c(x)$ for some $c(x)\in{\bf K}[x]$. Then $h_{(d-1)x}=c'(x)\in{\bf
K}[x]$. It follows from equation $(2.43)$ that $u_y=0$ or
$v_0'(u)\in {\bf
K}$.\\
If $u_y=0$, then it reduces to (I).\\
If $v_0'(u)\in {\bf K}$, then $v_0(u)=au+b$ for some $a,b\in {\bf
K}$. Thus, we have $v(u,h)=v_1h+au+b$ with $v_1\in {\bf K}^*$ and
$a,b\in {\bf K}$. Since $H(0)=0$, we have $b=0$. That is,
$v(u,h)=v_1h+au$. Thus, $u,v,h$ are linearly dependent.
\end{proof}

\begin{cor}
Let $H=(u(v,h),v(x,y),h(x,y,z))$ be a polynomial map with $H(0)=0$.
If $JH$ is nilpotent, then $u,v,h$ are linearly dependent.
\end{cor}
\begin{proof}
Let \[T=\left(
  \begin{array}{ccc}
    0 & 1 & 0 \\
    1 & 0 & 0 \\
    0 & 0 & 1 \\
  \end{array}
\right).\]  Then $T^{-1}HT=(v(y,x),u(v,h),h(y,x,z))$. Since $JH$ is
nilpotent, so $J(T^{-1}HT)=T^{-1}JHT$ is nilpotent. It follows from
Theorem 2.6 that $u,v,h$ are linearly dependent.
\end{proof}

\section{Polynomial maps of the form $H=(u(x,y),v(x,\allowbreak y,z), h(x,y))$}

In this section, we classify polynomial maps of the form
$H=(u(x,y),v(x,y,z), h(x,\allowbreak y))$ in the case $JH$ is
nilpotent and $(\deg_yu,\deg_yh)\leq 2$ or at least one of
$\deg_yu$, $\deg_yh$ is a prime.

\begin{lem}
Let $u(x,y),h(x,y)\in {\bf K}[x,y]$ and $u(0,0)=h(0,0)=0$. If $\det
J(u,h)=0$, then there exists $q(x,y)\in {\bf K}[x,y]$ such that
$u(x,y), h(x,y)\in {\bf K}[q(x,y)]$.
\end{lem}
\begin{proof}
Let $D=h_y\partial_x-h_x\partial_y$. It follows from Theorem 2.8 in
\cite{11} or Theorem 1.2.5 in \cite{3} that $\operatorname{Ker}
D={\bf K}[q]$ for some polynomial $q\in {\bf K}[x,y]$. Since
$Dh=0=Du$, so we have that $u,h\in {\bf K}[q]$.
\end{proof}

\begin{lem}
Let $H=(u(x,y),v(x,y,z),h(x,y))$ be a polynomial map with $H(0)=0$.
Assume that the components of $H$ are linearly independent over
${\bf K}$. If $JH$ is nilpotent, then $\deg_zv(x,y,z)=1$ and the
coefficient of $z$ in $v$ is a non-zero constant.
\end{lem}
\begin{proof}
The conclusion follows from the first part of the proof of Theorem
2.8 in \cite{13}.
\end{proof}

\begin{thm}
Let $H=(u(x,y),v(x,y,z),h(x,y))$ be a polynomial map with $H(0)=0$.
Assume that the components of $H$ are linearly independent over
${\bf K}$. If $JH$ is nilpotent and $\deg_yu$ is a prime or
$\deg_yh$ is a prime, then $u=g(ay+b(x))$,
$v=v_1z-a^{-1}b'(x)g(ay+b(x))-v_1l_2x$, $h=c_0u^2+l_2u$, where
$b(x)=v_1c_0ax^2+l_1x+\tilde{l}_2$; $v_1,c_0,a\in {\bf K}^*$;
$l_1,l_2,\tilde{l}_2\in {\bf K}$, $g(t)\in {\bf K}[t]$ and $g(0)=0$,
$\deg_tg(t)\geq 1$.
\end{thm}
\begin{proof}
Since $JH$ is nilpotent, we have the following equations:
\begin{equation}
\nonumber
  \left\{ \begin{aligned}
  u_x+v_y = 0~~~~~~~~~~~~~~~~~~(3.1) \\
  u_xv_y-v_xu_y = v_zh_y~~~~~~~~(3.2) \\
  v_z(u_xh_y-u_yh_x) = 0~~~~~~~(3.3)
                          \end{aligned} \right.
  \end{equation}
Let $v=v_dz^d+\cdots+v_1z+v_0$. Then it follows from Lemma 3.1 that
$d=1$ and $v_1\in {\bf K}^*$. Since $\deg_yu$ or $\deg_yh$ is a
prime, we have that $(\deg_yu,\deg_yh)=1$ or $\deg_yu$ or $\deg_yh$.

Case (I) If $(\deg_yu,\deg_yh)=1$, then the conclusion follows from
Theorem 2.8 in \cite{13}.

Case (II) If $(\deg_yu, \deg_yh)=\deg_yu$, then it follows from
equation $(3.3)$ and Lemma 3.1 that there exists $q\in {\bf K}[x,y]$
such that $u ,h\in {\bf K}[q]$. Since $\deg_yu$ is a prime, we have
that $\deg_yq=1$ or $\deg_yu$.\\
If $\deg_yq=1$, then the conclusion follows from the proof of
Theorem 2.8 in \cite{13}.\\
If $\deg_yq=\deg_yu$, then $u(x,y)=u(q)=\lambda q+\lambda_0$ for
$\lambda \in {\bf K}^*$, $\lambda_0\in {\bf K}$. Thus,
$q=\lambda^{-1}u-\lambda^{-1}\lambda_0$. That is, $h$ is a
polynomial of $u$. Then the conclusion follows from Theorem 2.1 in
\cite{18}.

Case (III) If $(\deg_yu, \deg_yh)=\deg_yh$, then $\deg_yh$ is a
prime. Thus, it follows from the arguments of Case (II) that $u$ is
a polynomial of $h$. It follows from Corollary 2.3 that $u,v,h$ are
linearly dependent. This is a contradiction!
\end{proof}

\begin{lem}
Let $H=(u(x,y),v(x,y,z),h(x,y))$ be a polynomial map over ${\bf
K}[x,y,z]$. Assume that $H(0)=0$ and the components of $H$ are
linearly independent over ${\bf K}$. If $JH$ is nilpotent, then
$\deg h\leq 2\deg u$, $\deg_yh\leq 2\deg_yu$, $\deg_xh\leq 2\deg_xu$
and the coefficients of the highest degree of $y$ in $u$ and $h$ are
non-zero constants.
\end{lem}
\begin{proof}
Let $v=v_dz^d+v_{d-1}z^{d-1}+\cdots+v_1z+v_0$. Then it follows from
Lemma 3.2 that $d=1$ and $v_1\in {\bf K}^*$. Since $JH$ is
nilpotent, so we have the following equations:
\begin{equation}
\nonumber
  \left\{ \begin{aligned}
  u_x+v_{0y} = 0~~~~~~~~~~~~~~~~~~~~~~~~(3.4) \\
  u_xv_{0y}-v_{0x}u_y-v_1h_y=0~~~~~~~~(3.5) \\
  v_1(u_xh_y-u_yh_x) = 0~~~~~~~~~~~~~~(3.6)
                          \end{aligned} \right.
  \end{equation}
It follows from equation $(3.4)$ that $u_x=-v_{0y}$. Thus, there
exists $P\in{\bf K}[x,y]$ such that
$$u=-P_y,~v_0=P_x~~~~~~~~~~~~~~~~~~~(3.7)$$
It follows from equation $(3.6)$ and Lemma 3.1 that there exists
$q\in {\bf K}[x,y]$ such that $u,h\in {\bf
K}[q]~~~~~~~~~~~~~~~~~~~~~~~~~~~~~~~~~~~~~~~~~~~~~~~~~~(3.8)$\\
It follows from $(3.7)$ and $(3.8)$ that $u_y=u'(q)q_y=-P_{yy}$,
$h_y=h'(q)q_y$, so we have
$$q_y=-\frac{P_{yy}}{u'(q)}~~~~~~~~~~~~~~~~~~~~~~~~~~~~(3.9)$$
and
$$h_y=-\frac{h'(q)}{u'(q)}P_{yy}~~~~~~~~~~~~~~~~~~~~~~(3.10)$$
because $u'(q)h'(q)\neq 0$. Otherwise, $u=0$ or $h=0$ which deduce
that $u,v,h$ are linearly dependent. This is a contradiction!
Substituting equations $(3.7)$, $(3.9)$ and $(3.10)$ to equation
$(3.5)$, we have the following equation:
$$u'(q)(P_{xy}^2-P_{xx}P_{yy})=v_1h'(q)P_{yy}~~~~~~~~~~(3.11)$$
Since $u,v,h$ are linearly independent, so $u_y\neq 0$ because $JH$
is nilpotent and $u(0,0)=0$. Thus, $P_{yy}=-u_y\neq 0$. Therefore,
we have the following inequality
$(\deg_qh(q)-1)\deg_yq+\deg_yP-2\leq
(\deg_qu(q)-1)\deg_yq+2(\deg_yP-1)$ by comparing the degree of $y$
of equation $(3.11)$. That is, $$\deg_qh(q)\deg_yq\leq
\deg_qu(q)\deg_yq+\deg_yP.$$ It follows from $(3.7)$ and $(3.8)$
that $\deg_yP-1=\deg_yu=\deg_qu(q)\deg_yq$. Thus, we have the
following inequality
$$\deg_yq(\deg_qh(q)-2\deg_qu(q))\leq 1~~~~~~~~~~~~~~~(3.12)$$
Since $u_y\neq 0$, so it follows from
$(3.8)$ that $\deg_yq\geq 1$.

If $\deg_yq=1$, then it follows from the proof of Theorem 2.8 in
\cite{13} that $h=c_0u^2+c_1u$ for $c_0\in {\bf K}^*$ and $c_1\in
{\bf K}$ and $u=g(ay+b(x))$ for $g(t)\in {\bf K}[t]$, $a\in {\bf
K}^*$.Then the conclusion follows.

If $\deg_yq\geq 2$, then it follows from $(3.12)$ that
$\deg_qh(q)\leq 2\deg_qu(q)$. Thus, we have $\deg h=\deg_qh(q)\deg
q\leq 2\deg_qu(q)\deg q=2\deg u$, $\deg_y h=\deg_qh(q)\deg_y q
\allowbreak \leq 2\deg_qu(q)\deg_y q=2\deg_y u$, $\deg_x
h=\deg_qh(q)\deg_x q\leq 2\deg_qu(q)\deg_x q=2\deg_x u$. Let
$P(x,y)=a_r(x)y^r+a_{r-1}(x)y^{r-1}+\cdots+a_1(x)y+a_0(x)$,
$h(x,y)=h_n(x)y^n+h_{n-1}(x)y^{n-1}+\cdots+h_1(x)y+h_0(x)$ with
$a_r(x)h_n(x)\neq 0$. It follows from equations $(3.5)$ and $(3.7)$
that\\
$(ra_r'(x)y^{r-1}+(r-1)a_{r-1}'(x)y^{r-2}+\cdots+a_1'(x))^2-(a_r''(x)y^r+a_{r-1}''(x)y^{r-1}+\cdots+a_1''(x)y+a_0''(x))
(r(r-1)a_r(x)y^{r-2}+(r-1)(r-2)a_{r-1}(x)y^{r-2}+\cdots+2a_2(x))=-v_1(nh_n(x)y^{n-1}+(n-1)h_{n-1}(x)y^{n-2}+\cdots+h_1(x))$
~~~~~~~~~~~~~~~(3.13)\\
We view that the polynomials are in ${\bf K}[x][y]$ with
coefficients in  ${\bf K}[x]$ when comparing the coefficients of
$y^k$. Since $n=\deg_yh\leq 2\deg_yu=2(r-1)$. So we have $n-1\leq
2r-3$. Then we have the following equation
$$r^2(a_r'(x))^2-r(r-1)a_r''(x)a_r(x)=0~~~~~~~~~~~~~~~(3.14)$$
by comparing the coefficients of $y^{2r-2}$ of equation $(3.13)$.
Thus, we have $a_r(x)\in {\bf K}^*$ by comparing the coefficients of
the highest degree of equation $(3.14)$. Let
$q(x,y)=q_l(x)y^l+\cdots+q_1(x)y+q_0(x)$. Then
$u(x,y)=u(q)=-P_y=-(ra_r(x)y^{r-1}+\cdots+a_1(x))$. Thus, we have
$q_l(x)\in {\bf K}^*$. Since
$h=h(q)=h(q_l(x)y^l+\cdots+q_1(x)y+q_0(x))$, so the coefficients of
the highest degree of $y$ in $u$ and $h$ are non-zero constants.
\end{proof}

\begin{thm}
Let $H=(u(x,y),v(x,y,z),h(x,y))$ be a polynomial map over ${\bf
K}[x,y,z]$. Assume that $H(0)=0$ and the components of $H$ are
linearly independent over ${\bf K}$. If $JH$ is nilpotent and
$(\deg_yu, \deg_yh)\leq 2$, then $H$ has the form of Theorem 3.3.
\end{thm}

\begin{proof}
Let $v=v_dz^d+\cdots+v_1z+v_0$. Then it follows from Lemma 3.2 that
$d=1$ and $v_1\in{\bf K}^*$. Since $JH$ is nilpotent, we have the
following equations:
\begin{equation}
\nonumber
  \left\{ \begin{aligned}
  u_x+v_{0y} = 0~~~~~~~~~~~~~~~~~~~~~~~~(3.4) \\
  u_xv_{0y}-v_{0x}u_y-v_1h_y=0~~~~~~~~(3.5) \\
  v_1(u_xh_y-u_yh_x) = 0~~~~~~~~~~~~~~(3.6)
                          \end{aligned} \right.
  \end{equation}
It follows from equation $(3.6)$ and Lemma 3.1 that there exists
$q\in {\bf K}[x,y]$ such that $u,h\in {\bf
K}[q]~~~~~~~~~~~~~~~~~~~~~~~~~~~~~~~~~~~~~~~~~~~~~~~~~~(3.8)$\\
Since $\deg_yq\leq (\deg_yu,\deg_yh)\leq 2$, so we have $\deg_yq=0$
or 1 or 2. \\
If $\deg_yq=0$, then $\deg_yu=0=\deg_yh$. It follows from equation
$(3.5)$ that $v_{0y}=0$ or $u_x=0$. If $v_{0y}=0$, then it follows
from equation $(3.4)$ that $u_x=0$. Since $u(0,0)=0$, we have $u=0$
in the two cases. Thus, $u,v,h$ are linearly dependent. This is a
contradiction!\\
If $\deg_yq=1$, then the conclusion follows from the proof of
Theorem 2.8 in \cite{13}.\\
If $\deg_yq=2$, then $\deg_yq_y=1$. Let
$q(x,y)=q_2(x)y^2+q_1(x)y+q_0(x)$ with $q_2(x)\neq 0$. It follows
from Lemma 3.4 and $(3.8)$ that $q_2(x)\in {\bf K}^*$. Thus,
$q_y=2q_2y+q_1(x)$ with $q_2\in {\bf K}^*$. Clearly, $q_y$ is
irreducible.

Substituting $(3.4)$ and $(3.8)$ to equation $(3.5)$, we have
$-(u'(q)q_x)^2-v_{0x}u'(q)q_y=v_1h'(q)q_y$. That is,
$$q_y[v_1h'(q)+v_{0x}u'(q)]=-(u'(q)q_x)^2$$
Since $q_y$ is irreducible, so we  have that $q_y|q_x$ or
$q_y|u'(q)$.

If $q_y|q_x$, then there exists $\mu(x,y)\in {\bf K}[x,y]$ such that
$q_x=\mu(x,y)q_y$. Since $q_y=2q_2y+q_1(x)$ and
$q_x=q_1'(x)y+q_0'(x)$, so we have $\mu(x,y)\in {\bf K}[x]$. Let
$\mu(x):=\mu(x,0)=\mu(x,y)$. Then
$$q_x=\mu(x)q_y~~~~~~~~~~~~~~~~~~~(3.15)$$
Let $\bar{x}=x$, $\bar{y}=y+\int\mu(x)dx$. Then it follows from
equation $(3.15)$ that $q_{\bar{x}}=0$. Let
$\bar{b}(x)=\int\mu(x)dx$. Then $q=q(y+\bar{b}(x))\in {\bf
K}[y+\bar{b}(x)]$. That is, $u,h\in {\bf K}[y+\bar{b}(x)]$. Then the
conclusion follows from the proof of Theorem 2.8 in \cite{13}.

If $q_y|u'(q)$, then $(2q_2y+q_1(x))|u'(q)$. Since $u'(q)$ is a
polynomial of $q$, so we have that
$u'(q)=c_0(q+c_1)(q+c_2)\cdots(q+c_k)$ for $c_0\in {\bf K}^*$,
$c_i\in{\bf K}$, $1\leq i\leq k$. Since $q_y$ is irreducible, so
there exists $i_0\in\{1,2,\ldots,k\}$ such that $q_y|(q+c_{i_0})$.
That is,
$$q_2y^2+q_1(x)y+q_0(x)+c_{i_0}=(cy+\alpha(x))(2q_2y+q_1(x))~~~~~~~~~~~(3.16)$$
Then we have the following equations:
\begin{equation}
\nonumber
  \left\{ \begin{aligned}
  c = \frac{1}{2}~~~~~~~~~~~~~~~~~~~~~~~~~~~~~~~~(3.17) \\
  cq_1(x)+2q_2\alpha(x)=q_1(x)~~~~~~~(3.18) \\
 \alpha(x)q_1(x) = q_0(x)+c_{i_0}~~~~~~~~~(3.19)
                          \end{aligned} \right.
  \end{equation}
by comparing the coefficients of $y^2$, $y$, $y^0$ of equation
$(3.16)$. It follows from equations $(3.18)$, $(3.19)$ that
$\alpha(x)=\frac{1}{4q_2}q_1(x)$,
$q_0(x)=\frac{1}{4q_2}q_1^2(x)-c_{i_0}$. Thus, we have
$q(x,y)=q_2y^2+q_1(x)y+\frac{1}{4q_2}q_1^2(x)-c_{i_0}=q_2\cdot
(y+\frac{1}{2q_2}q_1(x))^2-c_{i_0}$. That is, $q(x,y)\in {\bf
K}[y+\frac{q_1(x)}{2q_2}]$. Thus, we have that $u,h\in {\bf
K}[y+\frac{q_1(x)}{2q_2}]$. Then the conclusion follows from the
proof of Theorem 2.8 in \cite{13}.
\end{proof}

\begin{rem}
We can replace the condition that $(\deg_yu(x,y),m)=1$ by the
condition $(\deg_yu(x,y),m)\leq 2$ in Theorem 2.10 and replace the
condition that $(m,n)=1$ by the condition $(m,n)\leq 2$ in Theorem
3.2 and Theorem 3.4 in \cite{13}.
\end{rem}

\begin{cor}
Let $H=(u(x,y,z),v(x,y),h(x,y))$ be a polynomial map over ${\bf
K}[x,y,z]$. Assume that $H(0)=0$ and the components of $H$ are
linearly independent over ${\bf K}$. If $JH$ is nilpotent and
$(\deg_xv, \deg_xh)\leq 2$ or at least one of $\deg_xv$, $\deg_xh$
is a prime, then $u=u_1z-a^{-1}b'(y)g(ax+b(y))-u_1l_2y$,
$v=g(ax+b(y))$, $h=c_0v^2+l_2v$, where
$b(y)=u_1c_0ay^2+l_1y+\tilde{l}_2$, $u_1, c_0, a\in{\bf K}^*$,
$l_1,l_2, \tilde{l}_2\in{\bf K}$, $g(t)\in{\bf K}[t]$ and $g(0)=0$,
$\deg_tg(t)\geq 1$.
\end{cor}
\begin{proof}
Let \[T=\left(
  \begin{array}{ccc}
    0 & 1 & 0 \\
    1 & 0 & 0 \\
    0 & 0 & 1 \\
  \end{array}
\right).\] Then
$T^{-1}HT=(v(y,x),u(y,x,z),h(y,x)):=(\bar{v}(x,y),\bar{u}(x,y,z),\bar{h}(x,y))$.
Since $JH$ is nilpotent, so we have that $J(T^{-1}HT)=T^{-1}JHT$ is
nilpotent. Clearly, $\deg_y\bar{v}(x,y)=\deg_yv(y,x)=\deg_xv(x,y)$
and $\deg_y\bar{h}(x,y)=\deg_yh(y,x)=\deg_xh(x,y)$, so we have that
$(\deg_y\bar{v}(x,y),\deg_y\bar{h}(x,y))\leq 2$ or at least one of
$\deg_y\bar{v}$, $\deg_y\bar{h}$ is a prime. It follow from Theorem
3.3 and Theorem 3.5 that $v(y,x)=g(ay+b(x))$,
$u(y,x,z)=u_1z-a^{-1}b'(x)g(ay+b(x))-u_1l_2x$,
$h(y,x)=c_0v^2(y,x)+l_2v(y,x)$, where
$b(x)=u_1c_0ax^2+l_1x+\tilde{l}_2$; $u_1,c_0,a\in {\bf K}^*$;
$l_1,l_2,\tilde{l}_2\in {\bf K}$, $g(t)\in {\bf K}[t]$ and $g(0)=0$
$\deg_tg(t)\geq 1$. Then the conclusion follows.
\end{proof}

\begin{cor}
Let $H=(u(x,y,z),v(x,y),h(x,y))$ be a polynomial map over ${\bf
K}[x,y,z]$ and $v(x,y)=(a(y)x+b(y))^n$. Assume that $H(0)=0$ and the
components of $H$ are linearly independent over ${\bf K}$. If $JH$
is nilpotent, then $v(x,y)=(ay+b(y))^n$,
$u=u_1z-a^{-1}b'(y)(ax+b(y))^n-u_1l_2y$, $h=c_0v^2+l_2v$, where
$b(y)=u_1c_0ay^2+l_1y+\tilde{l}_2$, $u_1, c_0, a\in{\bf K}^*$,
$l_1,l_2, \tilde{l}_2\in{\bf K}$.
\end{cor}
\begin{proof}
Let  \[T=\left(
  \begin{array}{ccc}
    0 & 1 & 0 \\
    1 & 0 & 0 \\
    0 & 0 & 1 \\
  \end{array}
\right).\] Then
$T^{-1}HT=(v(y,x),u(y,x,z),h(y,x)):=(\bar{v}(x,y),\bar{u}(x,y,z),\bar{h}(x,y))$.
Since $JH$ is nilpotent, so we have that $J(T^{-1}HT)=T^{-1}JHT$ is
nilpotent. It follows from Lemma 3.2 that $\deg_z\bar{u}(x,y,z)=1$.
Thus, we have that $\bar{v}_x\bar{h}_y-\bar{v}_y\bar{h}_x=0$. It
follows from Lemma 3.1 that there exists $q\in{\bf K}[x,y]$ such
that $\bar{v},\bar{h}\in {\bf K}[q]$. Since
$\bar{v}(x,y)=v(y,x)=(a(x)y+b(x))^n$, we have $q=(a(x)y+b(x))^{n_1}$
for some $n_1\in {\bf N}^*$. Thus, we have $\bar{h}\in {\bf
K}[a(x)y+b(x)]$. That is, $h\in {\bf K}[a(y)x+b(y)]$. Since the
condition $(\deg_y\bar{v},\deg_y\bar{h})\leq 2$ in Theorem 3.5 is
only used to get that $\bar{v},\bar{h}$ are polynomials of
$a(x)y+b(x)$ and the condition $(\deg_xv,\deg_xh)\leq 2$ in
Corollary 3.7 is only used to get that
$(\deg_y\bar{v},\deg_y\bar{h})\leq 2$, so the conclusion follows
from the the proof the Theorem 3.5 and Corollary 3.7.
\end{proof}

\begin{cor}
Let $H=(u(x,y),v(x,y,z),h(x,y))$ be a polynomial map over ${\bf
K}[x,y,z]$. Assume that $H(0)=0$ and the components of $H$ are
linearly independent over ${\bf K}$. If $JH$ is nilpotent and
$\deg_yu\leq 4$ or $\deg_yh\leq 4$, then $H$ has the form of Theorem
3.3.
\end{cor}

\begin{proof}
Let $v=v_dz^d+\cdots+v_1z+v_0$. Then it follows from Lemma 3.2 that
$d=1$ and $v_1\in {\bf K}^*$. Since $JH$ is nilpotent, we have the
following equations:
\begin{equation}
\nonumber
  \left\{ \begin{aligned}
  u_x+v_{0y} = 0~~~~~~~~~~~~~~~~~~~~~~~~(3.4) \\
  u_xv_{0y}-v_{0x}u_y-v_1h_y=0~~~~~~~~(3.5) \\
  v_1(u_xh_y-u_yh_x) = 0~~~~~~~~~~~~~~(3.6)
                          \end{aligned} \right.
  \end{equation}
If $\deg_yu=3$ or $\deg_yh=3$, then the conclusion follows from
Theorem 3.3. So we can assume that $\deg_yu\neq 3$ and $\deg_yh\neq
3$. It follows from equation $(3.6)$ and Lemma 3.1 that there exist
$q\in{\bf K}[x,y]$ such that $u,h\in{\bf K}[q]$. Since
$\deg_yq|(\deg_yu,\deg_yh)$, so we have $\deg_yq=0$ or 1 or 2 or 4.

If  $\deg_yq=0$ or 1 or 2, then the conclusion follows from the
proof of Theorem 3.5.

If $\deg_yq=4$, then $\deg_yq=\deg_yu$ or $\deg_yq=\deg_yh$.

(1) If $\deg_yq=\deg_yu$, then $u(x,y)=u(q)=\lambda q+\lambda_0$
with $\lambda\in{\bf K}^*$, $\lambda_0\in {\bf K}$. That is,
$q=\lambda^{-1}u-\lambda^{-1}\lambda_0$. Thus, $h$ is a polynomial
of $u$. Then the conclusion follows from Theorem 2.1 in \cite{18}.

(2) If $\deg_yq=\deg_yh$, then it follows from the arguments of
$(1)$ that $u$ is a polynomial of $h$. It follows from Corollary 2.3
that $u,v,h$ are linearly dependent. This is a contradiction!
\end{proof}

\begin{cor}
Let $H=(u(x,y),v(x,y,z),h(x,y))$ be a polynomial map over ${\bf
K}[x,y,z]$. Assume that $H(0)=0$ and the components of $H$ are
linearly independent over ${\bf K}$. If $JH$ is nilpotent and $u$ or
$h$ is a polynomial of $y+a(x)$ for some $a(x)\in{\bf K}[x]$, then
$H$ has the form of Theorem 3.3.
\end{cor}
\begin{proof}
Let $v=v_dz^d+\cdots+v_1z+v_0$ with $v_i\in {\bf K}[x,y]$, $0\leq
i\leq d$. Then it follows from Lemma 3.2 that $d=1$ and $v_1\in {\bf
K}^*$. Since $JH$ is nilpotent, so we have the following equations:
\begin{equation}
\nonumber
  \left\{ \begin{aligned}
  u_x+v_{0y} = 0~~~~~~~~~~~~~~~~~~~~~~~~(3.4) \\
  u_xv_{0y}-v_{0x}u_y-v_1h_y=0~~~~~~~~(3.5) \\
  v_1(u_xh_y-u_yh_x) = 0~~~~~~~~~~~~~~(3.6)
                          \end{aligned} \right.
  \end{equation}
It follows from equation $(3.6)$ and Lemma 3.1 that there exists
$q\in {\bf K}[x,y]$ such that $u,h\in {\bf
K}[q]$~~~~~~~~~~~~~~~~~~~~~~~~~~~~~~~~~~~~~~~~~~~~~~~~~~~(3.8)\\
Thus, by the Fundamental Theorem of Algebra, we have that
$$u=u(q)=c_0(q+c_1)(q+c_2)\cdots(q+c_k)~~~~~~~(3.20)$$
and
$$h=h(q)=d_0(q+d_1)(q+d_2)\cdots(q+d_l)~~~~~~~~~~$$
for $c_0,d_0\in {\bf K}^*$ and $c_i,d_j\in {\bf K}$, $1\leq i\leq
k$, $1\leq j\leq l$.

Case I If $u$ is a polynomial of $y+a(x)$, then let $T=y+a(x)$, by
the Fundamental Theorem of Algebra, we have
$$u=u(T)=e_0(T+e_1)(T+e_2)\cdots(T+e_s)~~~~~~~~~~~(3.21)$$
for $e_0\in {\bf K}^*$, $e_1,\ldots,e_s\in {\bf K}$. It is clear
that $T+e_1$, $T+e_2$, $\ldots$ , $T+e_s$ are irreducible. It
follows from equations $(3.20)$ and $(3.21)$ that
$$(q+c_1)(q+c_2)\cdots(q+c_k)|(T+e_1)(T+e_2)\cdots(T+e_s)~~~~~~~~~(3.22)$$
Suppose $q+c_1=q_1q_2\cdots q_{r_1}$ and $q_1,q_2,\ldots, q_{r_1}$
are irreducible. It follows from $(3.22)$ that $q_m|(T+e_{i_m})$ for
$1\leq m\leq r_1$, $1\leq i_m\leq s$. That is, $T+e_{i_m}=b_mq_m$
for $b_m\in {\bf K}^*$. Thus, we have
$q+c_1=b(T+e_{i_1})(T+e_{i_2})\cdots(T+e_{i_{r_1}})$, where
$b=(b_1b_2\cdots b_m)^{-1}$. Therefore, $q$ is a polynomial of $T$.
So we have $q,u,h\in {\bf K}[y+a(x)]$. Then the conclusion follows
from the proof of Theorem 2.8 in \cite{13}.

Case II If $h$ is a polynomial of $y+a(x)$, then let $T=y+a(x)$, it
follows from the arguments of Case I that $u,h\in{\bf K}[y+a(x)]$.
Then the conclusion follows from the proof of Theorem 2.8 in
\cite{13}.
\end{proof}

\begin{thm}
Let $H=(u(x,y),v(x,y,z),h(x,y))$ be a polynomial map over ${\bf
K}[x,y,z]$. Assume that $H(0)=0$. If $JH$ is nilpotent and $u$ is
homogeneous of degree $n$, then $u,v,h$ are linearly dependent.
\end{thm}
\begin{proof}
Let $v=v_dz^d+\cdots+v_1z+v_0$ with $v_i\in {\bf K}[x,y]$, $0\leq
i\leq d$. Then it follows from Lemma 3.2 that $d=1$ and $v_1\in {\bf
K}^*$. Since $JH$ is nilpotent, so we have the following equations:
\begin{equation}
\nonumber
  \left\{ \begin{aligned}
  u_x+v_{0y} = 0~~~~~~~~~~~~~~~~~~~~~~~~(3.4) \\
  u_xv_{0y}-v_{0x}u_y-v_1h_y=0~~~~~~~~(3.5) \\
  v_1(u_xh_y-u_yh_x) = 0~~~~~~~~~~~~~~(3.6)
                          \end{aligned} \right.
  \end{equation}
It follows from equation $(3.4)$ that $u_x=-v_{0y}$. Thus, there
exists a polynomial $P\in {\bf K}[x,y]$ such that $u=-P_y,~ v_0=P_x$~~~~~~~~~~~~~~~~~~~~~~~~~~~~~~(3.23)\\
Let $h=h^{(m)}+h^{(m-1)}+\cdots+h^{(1)}$, where $h^{(j)}$ is the
homogeneous part of degree $j$ of $h$ and $h^{(m)}\neq 0$. It
follows from equation $(3.6)$ and Lemma 3.1 that there exists $q\in
{\bf K}[x,y]$ such that $u=u(q),~h=h(q)\in {\bf K}[q]$. Suppose
$\deg q=l$, $\deg_qu=s$, $\deg_qh=t$. Then $$\deg u=sl=n,~ \deg
h=tl=m~~~~~~~~~~~~~~~(3.24)$$ Let
$u=a_nx^n+a_{n-1}x^{n-1}y+a_{n-2}x^{n-2}y^2+\cdots+a_1xy^{n-1}+a_0y^n$.
It follows from equation $(3.23)$ that $P_y=-u$. Thus, we have
$$P=-P^{(n+1)}-f(x)~~~~~~~~~~~~~~~~~~(3.25)$$
for some $f(x)\in {\bf K}[x]$ and
$P^{(n+1)}=a_nx^ny+\frac{1}{2}a_{n-1}x^{n-1}y^2+\cdots+\frac{1}{n}a_1xy^{n}+\frac{1}{n+1}a_0y^{n+1}$.
Substituting equations $(3.23)$ and $(3.25)$ to equation $(3.5)$, we
have the following equation
$$-u_x^2+(P_{xx}^{(n+1)}+f''(x))u_y=v_1(h_y^{(m)}+h_y^{(m-1)}+\cdots+h_y^{(1)})~~~~~~(3.26)$$
Let $f''(x)=b_rx^r+b_{r-1}x^{r-1}+\cdots+b_1x+b_0$ with
$b_0,b_1,\ldots,b_r\in {\bf K}$ and $b_r\neq 0$.\\

(1) If $r\geq n$, then $m=n+r$ and $v_1h_y^{(m)}=b_rx^ru_y$. That
is, $$h^{(m)}=v_1^{-1}b_rx^ru+c(x)~~~~~~~~~~~~(3.27)$$ for some
$c(x)\in {\bf K}[x]$. The highest degree term in equation $(3.6)$ is
$u_xh_y^{(m)}-u_yh_x^{(m)}$. Thus, we have
$$u_xh_y^{(m)}-u_yh_x^{(m)}=0~~~~~~~~~~~~~~~~~~(3.28)$$
Substituting equation $(3.27)$ to equation $(3.28)$, we have that
$u_x(v_1^{-1}b_rx^ru_y)-u_y(v_1^{-1}b_rx^ru_x+rv_1^{-1}b_rx^{r-1}u+c'(x))=0$.
That is, $$(rv_1^{-1}b_rx^{r-1}u+c'(x))u_y=0.$$ Since
$rv_1^{-1}b_r\neq 0$, so we have $u_y=0$. Then $u_x=0$ because $JH$
is nilpotent. Since $u(0,0)=0$, so $u=0$. Thus, $u,v,h$ are linearly
dependent.\\

(2) If $r\leq n-1$ and
$-u_x^2+P_{xx}^{(n+1)}u_y+b_{n-1}x^{n-1}u_y\neq 0$ in the case
$r=n-1$ or $-u_x^2+P_{xx}^{(n+1)}u_y \neq 0$ in the case $r<n-1$,
then $$m=2n-1.$$ Substituting equation $(3.24)$ to the above
equation, we have that $tl= 2sl-1$. That is, $(2s-t)l=1$. So we have
that $l=1$ and $$t=2s-1~~~~~~~~~~~~~~(3.29)$$ Since $\deg_yq\leq
\deg q=1$, so it follows from the proof of Theorem 2.8 in \cite{13}
that $H$ has the form of Theorem 3.3. Therefore, $\deg h=2\deg u$.
It follows from equation $(3.24)$ that $t=2s$. This contradicted
with equation $(3.29)$.\\

(3) If $r\leq n-1$ and
$-u_x^2+P_{xx}^{(n+1)}u_y+b_{n-1}x^{n-1}u_y=0$ in the case $r=n-1$
or $-u_x^2+P_{xx}^{(n+1)}u_y = 0$ in the case $r<n-1$, then equation
(3.26) has the following form:
$$\bar{f}(x)u_y=v_1h_y~~~~~~~~~~~~~~(3.30)$$
where $\bar{f}(x)\in{\bf K}[x]$ and $\bar{f}(x)=f''(x)$ in the case
$r<n-1$, $\bar{f}(x)=f''(x)-b_{n-1}x^{n-1}$ in the case $r=n-1$. It
follows from equation $(3.30)$ that $h=\bar{\bar{f}}(x)u+\bar{c}(x)$
for $\bar{c}(x)\in {\bf K}[x]$,
$\bar{\bar{f}}(x)=v_1^{-1}\bar{f}(x)$. It follows from equation
$(3.6)$ that
$u_x(\bar{\bar{f}}(x)u_y)-u_y(\bar{\bar{f}}'(x)u+\bar{\bar{f}}(x)u_x+\bar{c}'(x))=0$.
That is, $(\bar{\bar{f}}'(x)u+\bar{c}'(x))u_y=0$. Thus, we have that
$u_y=0$ or $\bar{\bar{f}}'(x)=\bar{c}'(x)=0$.

If $u_y=0$, then it follows from the arguments of (1) that $u,v,h$
are linearly dependent.

If $\bar{\bar{f}}'(x)=\bar{c}'(x)=0$, then $\bar{\bar{f}}(x),
\bar{c}(x)\in {\bf K}$. That is, $h=\bar{\bar{f}}u+\bar{c}$. Since
$u(0,0)=h(0,0)=0$, so we have $\bar{c}=0$. Thus, $u,v,h$ are
linearly dependent.
\end{proof}


\begin{thebibliography}{99}
\bibitem{2} H. Bass, E. Connell and D. Wright, \newblock {\em The Jacobian Conjecture: Reduction of Degree and Formal Expansion of the Inverse}, Bulletin of the American Mathematical Society,\@ 7 (1982), 287-330.
\bibitem{10} M. de Bondt, \newblock {\em Quasi-translations and counterexamples to the homogeneous dependence problem}, Proceedings of the American Mathematical Society 134 (2006) 2849-2856.
\bibitem{8} M. de Bondt and A. van den Essen, \newblock {\em The Jacobian conjecture: linear triangularization for homogeneous polynomial maps in dimension three}, Report 0413, University of Nijmegen, The Netherlands, 2004.
\bibitem{18} Marc Chamberland and Arno van den Essen, \newblock {\em Nilpotent Jacobian in dimension three}, Journal of Pure and Applied Algebra \@ 205 (2006) 146-155.
\bibitem{3} A. van den Essen, Polynomial Automophisms and the Jacobian Conjecture, Vol. \@ 190 in Progress in Mathematics Birkhauser Basel, \@ 2000.
\bibitem{9} A. van den Essen, \newblock {\em Nilpotent Jacobian matrices with independent rows}, Report 9603, University of Nijmegen, The Netherlands, 1996.
\bibitem{7} E. Hubbers, \newblock {\em The Jacobian conjecture: cubic homogeneous maps in dimension four}, Master's Thesis, University of Nijmegen, The Netherlands, 1994.
\bibitem{1} O.H. Keller, \newblock {\em Ganze Cremona-transformationen Monatschr.}, Math. Phys., \@ 47 (1939) pp.\@ 229-306.
\bibitem{11} A. Nowicki, M. Nagata, \newblock {\em Rings of constants for k-derivations in
$k[x_1,\ldots,x_n]$}, J. Math. Kyoto Univ. 28(1988), No. 1, 111-118.
\bibitem{4} S.S.S.Wang, \newblock {\em A Jacobian criterion for separability}, Jour. of Algebra,\@ 65 (1980), 453-494.
\bibitem{6} D. Wright, \newblock {\em The Jacobian conjecture: linear triangularization for cubics in dimension three}, Linear and Multilinear Algebra 34(1993) 85-97.
\bibitem{5} A.V.Yagzhev, \newblock {\em On Keller's problem}, Siberian Math. Journal, \@ 21 (1980), 747-754.
\bibitem{13} Dan Yan, Guoping Tang, \newblock {\em Polynomial maps with nilpotent Jacobians in dimension
three}, Linear Algebra and its Applications, \@ 489 (2016), 298-323.
\end{thebibliography}
\end{document}